\documentclass[12pt]{amsart}

\usepackage{amsmath,amscd,amssymb}

\usepackage{latexsym}
\usepackage{amsfonts}
\usepackage[all]{xy}

\begin{document}
\title{Non-Abelian Localization for $U(1)$ Chern-Simons Theory}
\date{January 5, 2009}

\author{Lisa Jeffrey}
\address{Department of Mathematics, University of Toronto, Canada}
\email{jeffrey@math.toronto.edu}

\author{Brendan McLellan}
\address{Department of Mathematics, University of Toronto, Canada}
\email{mclellan@math.toronto.edu}

\maketitle

\newtheorem{thm}{Theorem}[section]
\newtheorem{cor}[thm]{Corollary}
\newtheorem{lem}[thm]{Lemma}
\newtheorem{prop}[thm]{Proposition}
\newtheorem{define}{Definition}
\newtheorem{rmk}{Remark}
\newtheorem{example}{Example}
\newtheorem{conj}{Conjecture}
\newcommand{\dee}{\text{d}}
\newcommand{\Hom}{\mathrm{Hom}}
\newcommand{\Z}{\mathbb{Z}}
\newcommand{\R}{\mathbb{R}}
\newcommand{\C}{\mathbb{C}}
\newcommand{\I}{\mathbb{I}}
\newcommand{\Q}{\mathbb{Q}}
\newcommand{\U}{\mathbb{U}}
\newcommand{\F}{\mathbb{F}}
\newcommand{\X}{\overline{X}}
\newcommand{\gau}{\mathcal{G}}
\newcommand{\Tors}{\mathop{\mathrm{Tors}}\nolimits}









\section{Introduction}
\begin{center}
\textit{Dedicated to Hans Duistermaat on the occasion of his 65th birthday.}
\end{center}

This article studies the non-abelian localization results of
Beasley and Witten in  \cite{bw}, and considers the
analogue of these results when  the gauge group $G$ is the abelian
group $G=U(1)$.  In another case, Duistermaat and Heckman studied the stationary phase approximation for integrals of the form
\begin{equation*}
\int_{M}e^{it\mu^{Y}(x)}\beta(x)
\end{equation*}
where $\beta(x)=\frac{\omega^l}{l!(2\pi)^l}$, and $(M,\omega,G,\mu)$ is a Hamiltonian $G$-space, and $Y$ is an element of the
Lie algebra of $G$ and $\mu^Y$ the component of the moment map
in the direction of $Y$. They proved that the stationary phase approximation in this case is exact, and in the special case of isolated fixed points they obtained the (abelian) localization formula
\begin{equation*}{
\int_{M}e^{it\mu^{Y}(x)}\beta(x)=\left(\frac{i}{t}\right)^l\sum_{p\in M^{Y}}\frac{e^{it\mu^{Y}(p)}}{\sqrt{\text{det}\mathcal{L}_{p}(Y)}} }
\end{equation*}
where $M^{Y}=\{p\in M \,|\,Y_{M}(p)=0\}$ is the fixed point set of the vector field generated by $Y$.

  The case of $U(1)$ Chern-Simons theory is another situation where the stationary phase approximation is exact.  In \cite{bw}, Beasley and Witten study Chern-Simons gauge theory on a Seifert manifold $X$, with a gauge group $G$ that is non-abelian, compact, connected, simply connected and simple.
These assumptions imply that the principal $G$-bundle over $X$
\begin{displaymath}
\xymatrix{G \ar@{^{(}->}[r] & P \ar[d]\\
                              & X}
\end{displaymath}
is trivial.  This is not the case for $G=U(1)$, for which some of these
assumptions are not valid.

The authors of \cite{bw} then apply the technique of non-abelian localization to the Chern-Simons path integral,
\begin{equation}
\smash[t]{
Z_{X}(k)=\int_{\mathcal{M}}\mathcal{D}A\,\, \text{exp}\left[i\frac{k}{4\pi}\int_{X}\text{Tr}\left(A\wedge dA+\frac{2}{3}A\wedge A\wedge A\right)\right] }
\end{equation}
Here
 $\mathcal{M}$   is the moduli space of \textit{all} connections on $X$;
through localization, the authors of
\cite{bw} reformulate the partition function
as an integral over   the space of gauge equivalence
classes of  flat connections. They
are able to compute this partition function in terms of topological data on the moduli space of flat connections, $\mathcal{M}_0$, in several cases,
specifically related to    $SU(2)$.

Some results for the $U(1)$ Chern-Simons theory are already known --
see \cite{m} and \cite{mpr}.  We study
the results of Manoliu in \cite{m} where the partition function for $U(1)$
Chern-Simons theory has been calculated independently, and  compare her
 results to those of \cite{bw}.  Manoliu studies $U(1)$ Chern-Simons theory for
arbitrary 3-manifolds; on the other hand, Beasley and Witten study
Chern-Simons theory for simple simply connected gauge groups $G$ for
Seifert 3-manifolds (manifolds which are the total space of a  circle bundle
over a 2-manifold or a 2-dimensional orbifold). As noted above, one of the main differences of the $U(1)$ theory from the setting of \cite{bw} is the fact that there exist nontrivial principal $U(1)$-bundles over $X$.  This difference occurs explicitly in the expression for the partition function in the $U(1)$ case.
Manoliu gives the following formula for the
Chern-Simons partition function  in \cite{m}:

\begin{equation}
\smash[t]{ Z_{X} = \frac{k^{m_{X}}}{|\Tors H^{2}(X,\Z)|}\sum_{p\in \Tors H^{2}(X,\Z)}\sigma_{X,p}\int_{\mathcal{M}_0}(T_{X})^{1/2} }
\end{equation}
where 
$\text{Tors} H^{2}(X,\Z)$
is the torsion subgroup of $H^{2}(X,\Z)$ and is isomorphic to
$H^1(X,\R/\Z)$.
Note that a principal $U(1)$-bundle $P\rightarrow X$ has flat connections if and only if
the first Chern class $c_{1}(P)$ -- if non-zero -- is a torsion class in $H^2(X,\Z)$.
The sum over $p$ passes over different
topological types for the bundle over the 3-manifold  $X$.
We have calculated the dependence of the results of both
\cite{bw} and \cite{m} on the Chern-Simons parameter $k$.
  In our comparison we are only looking at the component of the partition function that comes from the contribution of
the trivial bundle  $P = U(1) \times X$  in the partition function of \cite{m}, since the results of \cite{bw} apply only to this case. For groups
$G$ which are simply connected, all $G$-bundles over three-manifolds
are trivial, whereas this is not true for $U(1)$ bundles.

We compare the expressions for the
Chern-Simons partition functions found by Beasley--Witten and
by Manoliu in their respective situations.
Both partition functions are expressed as integrals
of top degree differential forms over the
moduli space $\mathcal{M}$ of gauge
equivalence classes of flat connections over the 3-manifold.  We restrict to $G=U(1)$,
and compare
\begin{itemize}
\item the power of $k$ appearing in the integrand
\item the integrand (in both cases it is the symplectic volume form)
\end{itemize}
\begin{prop}
The $U(1)$ Chern-Simons partition function  from \cite{m} (Eq. 5.17)
and the specialization of Beasley and Witten's Chern-Simons partition
function (\cite{bw}, Eq. 5.172) to $G=U(1)$
are both proportional to
$$\smash[t]{k^{(2g-1)/2} \int_{U(1)^{2g} \times \Z_n} \omega^g }$$
(for the symplectic form $\omega$ on $U(1)^{2g}$).
\end{prop}

An important ingredient in our study of the the partition functions is the appearance of the Reidemeister torsion(R-torsion) in \cite{m}.  We provide a section devoted specifically to the R-torsion, and also study the relationship of the R-torsion to the symplectic volume of $\mathcal{M}_0$.


\section{$k$-dependence}

In this section we compare the results of \cite{bw} and \cite{m} for  the $k$-dependence of their partition functions.  In particular, we look at Eq. 5.172 of \cite{bw}, and Eq. 5.17 of \cite{m}, and derive their respective dependences on the Chern-Simons coupling constant $k$. 

Let us begin with \cite{m}.  Eq. 5.17 of \cite{m} reads:

\begin{align}
& Z_{X} = \frac{k^{m_{X}}}{|\Tors H^{2}(M,\Z)|} \sum_{p\in \Tors H^{2}(M,\Z) }\sigma_{X,p}\int_{\mathcal{M}_0} (T_{X})^{1/2} \\
      &= k^{m_{X}}\sum_{p\in \Tors H^{2}(M,\Z) }\int_{\mathcal{M}_{0,p}} \sigma_{X,p}(T_{X})^{1/2}\\
      &= k^{m_{X}}\int_{\mathcal{M}_{0}   }\sigma_{X}(T_{X})^{1/2}
\end{align}
where
\begin{equation} \label{e:exponent}
\smash[t]{m_{X}=\frac{1}{2}(\text{dim} H^{1}(X,\R)-\text{dim} H^{0}(X,\R))}
\end{equation}

Here if $p$ corresponds to a trivial
bundle $P = U(1) \times X $,   $\sigma_{X,p}(A) = e^{ikCS(A)}$ is the Chern-Simons
function of the connection $A$, exponentiated and to the power $k$.
We note that if $P$ is a trivial bundle and  $A$ is a critical point
for the Chern-Simons functional, then $dA = 0$ so
$\sigma_{X,p}(A) = 1$.

The $k$ dependence comes only
from the factor $k^{m_{X}}$. The value of  $m_{X}$ is as follows.
As $X$ is connected, $\text{dim} H^{0}(X,\R)=1$.
The values of $\text{dim} H^{1}(X,\R)$ are stated  in \cite{fs} - for completeness we provide a short proof.

\begin{prop}\cite{fs}
\begin{equation*}
{\rm dim} H^{1}(X,\R)=
\begin{cases}
2g & n\geq 1\\
2g+1 & n=0
\end{cases}
\end{equation*}
where $n$ is the degree of the $U(1)$-bundle $X$, and $g$ is the genus of the base space $\Sigma$
\begin{displaymath}
\xymatrix{U(1) \ar@{^{(}->}[r] & X \ar[d]\\
                              & \Sigma}
\end{displaymath}
\end{prop}

\begin{proof}
By the Universal Coefficient Theorem (UCT),
\begin{equation}\label{uct}
\smash[t]{H^{1}(X,\R)\simeq \text{Hom}(H_{1}(X,\Z),\R) }
\end{equation}
i.e. the UCT implies that:
\begin{equation*}
\smash[t]{0\rightarrow {\rm Ext} (H_{0}(X,\Z),\R)\rightarrow H^{1}(X,\R)\rightarrow \text{Hom}(H_{1}(X,\Z),\R)\rightarrow 0 }
\end{equation*}
is exact.  Also
\begin{equation*}
\smash[t]{{\rm Ext} (H_{0}(X,\Z),\R)\simeq {\rm Ext}(\Z,\R)\simeq 0}
\end{equation*}
since $\Z$ is free.

Thus we compute $\text{Hom}(H_{1}(X,\Z),\R)$.  By Hurewicz,
\begin{equation*}
\smash[t]{H_{1}(X,\Z)\simeq\frac{\pi_{1}(X)}{[\pi_{1}(X),\pi_{1}(X)]} }
\end{equation*}
We assume $X$ is a Seifert fibered manifold which fibers over a 2-manifold
 rather than over an orbifold, since this is the
setting of (5.172) of \cite{bw}.  
Hence we have the following presentation of $\pi_{1}(X)$(\cite{o1}):
\begin{equation*}
\smash[t]{ \pi_{1}(X)\simeq \langle a_{p},b_{p},h | [ a_{p},h]=[ b_{p},h]=1, \prod_{p=1}^{g}[ a_{p},b_{p}]=h^{n} \rangle }
\end{equation*}
where $g$ is the genus of the base space $\Sigma$ of our Seifert fibered 3-manifold $X$,
\begin{displaymath}
\xymatrix{U(1) \ar@{^{(}->}[r] & X \ar[d]\\
                              & \Sigma}
\end{displaymath}
and $n=c_{1}(X)$ is the Chern number of the $U(1)$-bundle $X$.  The generator $h$ arises from the generic fibre over $\Sigma$.  Observe
that in the abelianization of $\pi_{1}(X)$ the following relation
is satisfied:
\begin{equation}
\smash[t]{\prod_{p=1}^{g}[ a_{p},b_{p}]=h^{n}}
\end{equation}
We have
\begin{equation}
\smash[t]{[\pi_{1}(X),\pi_{1}(X)]=\langle [a_p,b_p]|\prod_{p=1}^{g}[ a_{p},b_{p}]=h^{n}\rangle
}\end{equation}
and therefore
\begin{equation}\smash[t]{
\frac{\pi_{1}(X)}{[\pi_{1}(X),\pi_{1}(X)]}=\langle a_p,b_p,h\,\,|\,\,[ a_{p},b_{p}]=h^{n}=1\rangle }
\end{equation} 
$$= \bigoplus_{p=1}^{g}\langle a_p \rangle \bigoplus_{p=1}^{g}\langle b_p \rangle \bigoplus \left(\frac{\langle h \rangle}{\langle h^{n} \rangle}\right) 
$$ 
where $a_p,b_p,h$ now represent equivalence classes in the abelianization and $\langle a_p \rangle\simeq \Z$, $\langle b_p \rangle\simeq \Z$, $\frac{\langle h \rangle}{\langle h^{n} \rangle}\simeq \Z/n\Z\simeq \Z_{n}$.
Thus,
\begin{equation}
\frac{\pi_{1}(X)}{[\pi_{1}(X),\pi_{1}(X)]}\simeq\begin{cases}
\Z^{2g}\times\Z_{n} & n\geq 1\\
\Z^{2g+1} & n=0\\
\end{cases}
\end{equation}
Finally we have
\begin{equation*} \smash[t]{
H^{1}(X,\R)\simeq \text{Hom}(H_{1}(X,\Z),\R)\simeq \text{Hom}\left(\frac{\pi_{1}(X)}{[\pi_{1}(X),\pi_{1}(X)]},\R\right) }
\end{equation*}
\begin{equation}
\simeq
\begin{cases}
\text{Hom}(\Z^{2g}\times\Z_{n}, \R) & n\geq 1\\
\text{Hom}(\Z^{2g+1}, \R) & n=0\\
\end{cases}
\end{equation}
\begin{equation}
\simeq
\begin{cases}
\R^{2g} & n\geq 1\\
\R^{2g+1} & n=0\\
\end{cases}
\end{equation}
In conclusion
\begin{equation}
\text{dim}H^{1}(X,\R)=
\begin{cases}
2g & n\geq 1\\
2g+1 & n=0
\end{cases}
\end{equation}
\end{proof}
We will restrict ourselves to the case $n \ne 0$
because this is assumed in \cite{bw}. Thus  we obtain the $k$ dependence
\begin{equation*}
\smash[t]{Z_{X}\sim k^{\frac{2g-1}{2}} }
\end{equation*}

Now consider the $k$ dependence for Eq. 5.172 of \cite{bw}.  For this computation we assume that all work leading to Eq. 5.172 is relevant to the case of a trivial principal $U(1)$-bundle over $X$.  We note that the main difference between our case  (with a $U(1)$ gauge group)
and the case studied in \cite{bw} is that \cite{bw} assume that their gauge group $G$ is compact, connected, simply connected and simple.  In particular, they can conclude that their principal $G$-bundle
\begin{displaymath}
\xymatrix{G \ar@{^{(}->}[r] & P \ar[d]\\
                              & X}
\end{displaymath}
is trivial.
 However $U(1)$ is neither simply connected or simple, and there exist non-trivial principal $U(1)$-bundles over $X$. It is surprising that our results show
\cite{bw} (5.172) is still valid in the $U(1)$ case although
we could not infer this from Beasley and Witten's calculation
because our situation does not satisfy the hypotheses of \cite{bw}.

Consider Eq. 5.172 of \cite{bw}:
\begin{equation}\label{bwf}
\smash[t]{Z_{X} := Z(\epsilon)|_{\widetilde{\mathcal{M}}_0} }
\end{equation}
$$\smash[t]{
=\frac{1}{|\Gamma|}
\text{exp}\left( -\frac{\imath\pi}{2}\eta_{0}\right)\int_{\widetilde{\mathcal{M}}_0}\widehat{A}(\widetilde{\mathcal{M}}_0)
\text{exp}\left[ \frac{1}{2\pi\epsilon}\Omega+\frac{1}{2}c_{1}(T\widetilde{\mathcal{M}}_0)+\frac{\imath n}{4\pi^{2}\epsilon_{r}}\Theta \right]
} $$
Here
\begin{itemize}
\item $\widetilde{\mathcal{M}}_0$ is a
smooth component of the moduli space of
irreducible flat connections on a Seifert manifold $X$ (we assume that our Seifert manifold $X$ is a smooth line bundle of degree $n$ over $\Sigma$)
\item $\Gamma=Z(G)$ is the center of $G$
\item $\eta_{0}=-\frac{n\text{dim} G}{6}$ 
\item $\widehat{A}(\widetilde{\mathcal{M}}_0)
=\prod^{\text{dim}\widetilde{\mathcal{M}}_0}_{j=1}\frac{x_{j}/2}
{\sinh(x_{j}/2)}$
 where $x_j   ( 1\leq j\leq n)$, are the Chern roots
of $T\widetilde{\mathcal{M}}_0$.
so that $c(T\widetilde{\mathcal{M}}_0)=\prod^{n}_{j=1}(1+x_{j}), \,\, x_{j}\in H^{2}(\widetilde{\mathcal{M}}_0,\Z)$
\item $\Omega$ is the symplectic form on $\widetilde{\mathcal{M}}_0$
\item $\epsilon_{r}=\frac{2\pi}{k+\widehat{c}_\mathfrak{g}}$, where $\widehat{c}_\mathfrak{g}$ is the dual Coxeter number of $G$
\item $\Theta\in H^{4}(\widetilde{\mathcal{M}}_0)$ is
the cohomology class corresponding to the degree 4 element $- (\phi, \phi)/2$
in the equivariant cohomology $H^4_G({\rm pt}) $
(for $\phi \in \mathfrak{g}$ using the Cartan model of equivariant
cohomology).
$\Theta$ can also be described in terms of the universal bundle $\U$
\begin{displaymath}
\xymatrix{\C \ar@{^{(}->}[r] & \U \ar[d]\\
                              & \text{Jac}(\Sigma)\times\Sigma}
\end{displaymath}
in other words
\begin{equation*}
\smash[t]{\Theta=-\frac{1}{2}c_{1}(\U)^{2}|_{\text{pt.}\in\Sigma}}
\end{equation*}
where $\text{Jac}(\Sigma)$ is the Jacobian of $\Sigma$.
\end{itemize}
The overall constant $\Gamma$ does not make sense for $G=U(1)$ since
it is infinite.
We
 disregard the overall constant  in front of the
integrand in $Z(\epsilon)|_{\widetilde{\mathcal{M}}_0}$.  Looking only at the $k$ dependence, we consider $Z(\epsilon)|_{\widetilde{\mathcal{M}}_0} $,
ignoring overall multiplicative constants:
\begin{equation}
\smash[t]{Z(\epsilon)|_{\widetilde{\mathcal{M}}_0}\sim \int_{\widetilde{\mathcal{M}}_0}\widehat{A}(\widetilde{\mathcal{M}}_0)\text{exp}\left[ \frac{1}{2\pi\epsilon}\Omega+\frac{1}{2}c_{1}(T\widetilde{\mathcal{M}}_0)+\frac{\imath n}{4\pi^{2}\epsilon_{r}}\Theta \right]
}\end{equation}
Note that $\widetilde{\mathcal{M}}_0\simeq U(1)^{2g}\times\Z_{n}$ by Proposition 2.2 of \cite{m}, and since we are restricting to the trivial bundle case, we identify $\widetilde{\mathcal{M}}_0\simeq U(1)^{2g}$ as the connected component corresponding to $p=0$.

The first thing we observe is that $\Theta=0$ in our case.  This follows since the universal bundle $\U$ for $U(1)$-bundles is the classical Poincare line bundle, and the Poincare line bundle is normalized to have degree $d=0$ when restricted to the Jacobian of $\Sigma$.  Since $c_{1}(\U)=d[\Sigma]\in H^{2}(\Sigma)$, this implies $c_{1}(\U)=0$, and hence $\Theta=0$.
Also, since $\widetilde{\mathcal{M}}_0\simeq U(1)^{2g}$, we know
\begin{align*}
\smash[t]{c(\widetilde{\mathcal{M}}_0)=c(T\widetilde{\mathcal{M}}_0)}
\end{align*}
                          $$\smash[t]{ =  \prod_{i=1}^{g}c(L_{i})=\prod_{i=1}^{g}(1+x_{i}) } $$
where $$\smash[t]{ L_{i}=T\Sigma_{i},\text{and}\,\, x_{i}=c_{1}(L_{i})\in H^2(\Sigma_{i}, \Z)
}$$

\noindent where $\Sigma_{i}\simeq(U(1))^2$.  Then the tangent bundles $T\Sigma_{i}$
are trivial, and hence
\begin{equation*}
\smash[t]{x_{i}=c_{1}(T\Sigma_{i})=0}
\end{equation*}
Thus
\begin{equation}
\smash[t]{\widehat{A}(\widetilde{\mathcal{M}}_0)=\prod^{\text{dim}\widetilde{\mathcal{M}}_0}_{j=1}\frac{x_{j}/2}{\text{sinh}(x_{j}/2)}=1}
\end{equation}
Clearly, $c_{1}(T\widetilde{\mathcal{M}}_0)=0$ as well, and we arrive at
\begin{equation}\label{symvol}
\smash[t]{Z(\epsilon)|_{\widetilde{\mathcal{M}}_0}\sim \int_{\widetilde{\mathcal{M}}_0}
\text{exp}\left[ \frac{1}{2\pi\epsilon}\Omega\right] }
\end{equation}
Recalling that $\epsilon=\frac{2\pi}{k}$, we have
\begin{align*}
\smash[t]{Z(\epsilon)|_{\widetilde{\mathcal{M}}_0} } &\sim \int_{\widetilde{\mathcal{M}}_0}
\text{exp}\left[ k\Omega\right]\\
                                      &=\int_{\widetilde{\mathcal{M}}_0}k^{g}\Omega^{g}/g!\\
                                      &=k^{g}\text{Vol}_{\Omega}(\widetilde{\mathcal{M}}_0)
\end{align*}
Thus the $U(1)$ Chern-Simons partition function
computed from Eq. 5.172 of \cite{bw} is
\begin{equation}
\smash[t]{Z_{X} := Z(\epsilon)|_{\widetilde{\mathcal{M}}_0}\sim k^{g}}
\end{equation}

There is a difference of
$k^{1/2}$ between the two cases, since
 $Z_{X} \sim k^{g}$ for Beasley-Witten, whereas $Z_{X} \sim k^{\frac{2g-1}{2}}$
for Manoliu.  Let us analyze this difference further.  In the case of \cite{m}, this extra factor of $k^{-1/2}$ appears because the dimension of the stabilizer of the gauge group action (for $U(1)$ gauge group) is $\text{dim}(H^0(X,\R))=1$.  

A similar phenomenon occurs in Yang-Mills theory at the higher 
non-flat critical points of the Yang-Mills action.
As observed in Section 4.3 of \cite{bw} (for example equation \cite{bw}(4.45)), there
 is a factor of $k^{1/2}$  in the Yang-Mills partition function coming from the fact that the gauge group 
$\mathcal{G}$ does not act locally freely 
on the locus of non-flat Yang-Mills solutions.  This $k^{1/2}$ factor
 comes from the $U(1)$ stabilizer at  a non-flat Yang-Mills solution.
  $U(1)$-Chern-Simons theory also has a $U(1)$ stabilizer at all points,  the subgroup of constant gauge transformations 
with values in $U(1)$.  
This accounts for the  extra factor of $k^{-1/2}$ in  the
Chern-Simons partition function in Manoliu's paper. 

In fact Beasley and Witten recast the Chern-Simons partition function
 as a Yang-Mills partition function (see (3.61) in \cite{bw}).
 In the computation of Eq. 5.172 of \cite{bw} it is assumed that one is 
localizing at an  irreducible flat connection, and therefore the isotropy group of $A$, $\Gamma_A=\{ u\in \mathcal{G}\,\,|\,\, u(A)=A\}$, is finite.
Hence there is no factor of  $k^{1/2}$ in \cite{bw} (5.172) because the dimension of the stabilizer is zero.  

The extra factor of $k^{-1/2}$ also appears in \cite{jkkw1} (see Section 9, Example 46).
This article treats integrals of the same
 form as \cite{bw} (3.61) over   symplectic manifolds equipped with Hamiltonian group actions. 
When the group acts locally freely at the zero locus of the moment map, it is shown in 
\cite{JK} 
that the integral is a polynomial in $\epsilon =  2\pi/k$. When the group
acts with nontrivial stabilizer at points in the zero locus of the moment map,
the partition function is a polynomial in $\sqrt{\epsilon}$ but not in $\epsilon$.


\section{Reidemeister Torsion and Symplectic Volume}
We would like to show that the remainder of the calculation in \cite{m} involving the Reidemeister torsion yields the symplectic volume as above.
In this section we review Reidemeister torsion (R-torsion) and provide some relevant examples.\\

The R-torsion is an invariant for a CW-complex and a representation of its fundamental group.  Before we define the R-torsion, we recall the definition of the torsion of a chain complex.
Let $$\smash[t]{C_{*}=\left(0\rightarrow C_{n}\xrightarrow{d_{n}} C_{n-1}\xrightarrow{d_{n-1}}\cdots C_{1}\xrightarrow{d_{1}} C_{0}\rightarrow 0\right)}$$
 be a chain complex over $\F$ (either $\R$ or $\C$).  Let $Z_{i}$ denote the cycles of this complex, $B_{i}$ denote the boundaries, and $H_{i}$ the homology.
Let $\{c^{i}\}$ be a basis of $C_{i}$ and $c$ be the collection $\{c^{i}\}_{i\geq 0}$.  We call a pair $(C_{*},c)$ a \textit{based chain complex}, $c$ the preferred basis of $C_{*}$ and $c^{i}$ the preferred basis of $C_{i}$.  Let $h^{i}$ be a basis of $H_{i}$.

We construct another basis as follows.  By the definitions of $Z_{i}$, $B_{i}$, and $H_{i}$, the following two split exact sequences exist:
\begin{eqnarray*}
&\smash[t]{ 0 \rightarrow Z_{i}\rightarrow C_{i}\xrightarrow{d_{i}} B_{i-1}\rightarrow 0,} &\\
& 0 \rightarrow B_{i}\rightarrow Z_{i}\rightarrow H_{i}\rightarrow 0.&
\end{eqnarray*}
Let $\widetilde{B}_{i-1}$ be a lift of $B_{i-1}$ to $C_{i}$ and $\widetilde{H}_{i}$ a lift of $H_{i}$ to $Z_{i}$.  Then we can decompose $C_{i}$ as follows.
\begin{align*}
C_{i} &= Z_{i}\oplus \widetilde{B}_{i-1}\\
      &= B_{i}\oplus \widetilde{H}_{i} \oplus \widetilde{B}_{i-1}\\
      &= d_{i+1} \widetilde{B}_{i}\oplus \widetilde{H}_{i}\oplus \widetilde{B}_{i-1}
\end{align*}
Choose a basis $b^{i}$ for $B_{i}$.  We write $\widetilde{b}^{i+1}=\{\widetilde{b}_{j}^{i+1}\}_{j=1}^{n_{i}}$ for a lift of $b^{i}$ and $\widetilde{h}^{i}=\{\widetilde{h}_{j}^{i}\}_{j=1}^{m_{i}}$ for a lift of $h^{i}$.  By construction, the set $\{\widetilde{b}^{i}\cup d_{i+1}(\widetilde{b}^{i+1}) \cup \widetilde{h}^{i}\}$ forms another ordered basis of $C_{i}$.  Denote this basis as $\{\widetilde{b}^{i} d_{i+1}(\widetilde{b}^{i+1})\widetilde{h}^{i}\}$.  The definition of the R-torsion, $\text{Tor}(C_{*},c,h)$, is as follows:
\begin{equation}\label{tor1}
\smash[t]{\text{Tor}(C_{*},c)\{h\}=(-1)^{|C_{*}|}\cdot \prod_{i=1}^{n}[\widetilde{b}^{i} d_{i+1}(\widetilde{b}^{i+1})\widetilde{h}^{i}/c^{i}]^{(-1)^{i+1}}\in\F^{*}
}\end{equation}
where $[\widetilde{b}^{i} d_{i+1}(\widetilde{b}^{i+1})\widetilde{h}^{i}/c^{i}]$ denotes the determinant of the change of basis matrix from the basis $\{c^{i}\}$
to the basis $\{\widetilde{b}^{i} d_{i+1}(\widetilde{b}^{i+1})\widetilde{h}^{i}\}$.
An alternative definition  is to equip our complex $C_{*}$ with volumes $\mu_{i}\in(\wedge^{\text{max}}C_{i})^{*}$, one for each $i$, and then to define
\begin{equation}\label{tor2}
\text{Tor}(C_{*},\mu)\{h\}=\frac{\bigwedge_{i \,\text{even}}\mu_{i}[\widetilde{b}^{i}\wedge d_{i+1}(\widetilde{b}^{i+1})\wedge\widetilde{h}^{i}]}{\bigwedge_{i \,\text{odd}}\mu_{i}[\widetilde{b}^{i}\wedge d_{i+1}(\widetilde{b}^{i+1})\wedge \widetilde{h}^{i}]}
\end{equation}
where we take $\widetilde{b}^{i}=\wedge_{j=1}^{n_{i}}\widetilde{b}^{i}_{j}$, $d_{i+1}\widetilde{b}^{i+1}=\wedge_{j=1}^{n_{i}}d_{i+1}\widetilde{b}^{i+1}_{j}$, and $\widetilde{h}^{i}=\wedge_{j=1}^{m_{i}}\widetilde{h}^{i}_{j}$.  The torsion is an element
\begin{equation*}
\smash[t]{\text{Tor}(C_{*},\mu)\in\otimes_{2i+1}[\wedge^{\text{max}}H_{2i+1}(C_{*})]\otimes_{2i}[\wedge^{\text{max}}H_{2i}(C_{*})]^{*} }
\end{equation*}
The latter definition specializes to the former definition when we choose the canonical volumes associated to a choice of preferred basis $c$ for $C$.

It is well known (see e.g. \cite{F}) that the torsion is independent of the
choices of $b^{i}$
and of the choices of  lifts $\widetilde{b^{i}}$, $\widetilde{h^{i}}$.
  In the case that our complex $C_{*}$ is acyclic, we define
\begin{equation}\label{tor3}
\smash[t]{\text{Tor}(C_{*},\mu)=\frac{\bigwedge_{i \,\text{even}}\mu_{i}[\widetilde{b}^{i}\wedge d_{i+1}(\widetilde{b}^{i+1})]}{\bigwedge_{i \,\text{odd}}\mu_{i}[\widetilde{b}^{i}\wedge d_{i+1}(\widetilde{b}^{i+1})]}\in \R^{*}. }
\end{equation}
We will be interested in a specific chain complex $C_*$.  In particular, let $N$ be a cell complex, and $\rho$ a representation of $\pi_{1}(N)$ in $G$.  The Lie algebra $\mathfrak{g}$ is acted on by $\pi_{1}(N)$ under the composition of the adjoint action of $G$ and the representation $\rho$. Let $\mathfrak{g}_{\rho}$ denote $\mathfrak{g}$ with the $\pi_{1}(N)$-module structure from $\rho$.  Let $\widetilde{N}$ denote the universal cover of $N$.  Since the fundamental group $\pi_{1}(N)$ acts on $\widetilde{N}$ by covering transformations, the chain complex $C_{*}(\widetilde{N})$ also has a natural $\pi_{1}(N)$-module structure.  The chain complex of interest is then $C_{*}(N,\mathfrak{g}_{\rho})$, defined as the quotient of $C_{*}(\widetilde{N})\otimes\mathfrak{g}$ under the equivalence:
\begin{equation}\label{eqrel}
\smash[t]{\sigma\otimes X \sim \sigma a\otimes \text{Ad}(\rho(a))^{-1}X}
\end{equation}
where $a\in \pi_{1}(N)$, $\sigma\in C_{*}(\widetilde{N})$, and $X\in\mathfrak{g}$.  The usual differential on $C_{*}(\widetilde{N})$ is compatible with the equivalence relation, and thus descends to a differential $\delta_{\rho}$ on $C_{*}(N,\mathfrak{g}_{\rho})$.  By dualizing one obtains the corresponding cochain complex $C^{*}(N,\mathfrak{g}_{\rho})$ with differential $d_{\rho}=\delta_{\rho}^{*}$.  We have the following

\begin{lem}\cite{jw}\label{conju}
Suppose $h\in G$.  If $\rho$ and $h\rho h^{-1}$ are conjugate representations of $\pi_{1}(N)$ in $G$, then the map $\text{Ad}(h):\mathfrak{g}\rightarrow\mathfrak{g}$ induces an isomorphism of the chain complexes $C_{*}(N,\mathfrak{g}_{\rho})$ and $C_{*}(N,\mathfrak{g}_{h\rho h^{-1}})$.  Hence one obtains a natural isomorphism between the cohomology groups $$H^{i}(C_{*}(N,\mathfrak{g}_{\rho}))$$ and $$H^{i}(C_{*}(N,\mathfrak{g}_{h\rho h^{-1}})).$$
\end{lem}

We will mainly be interested in the zeroth and first cohomology groups of this complex.  We recall the following
\begin{prop}\cite{jw}\label{jw}
Let $[\rho]\in\text{Hom}(\pi_{1}(N),G)/G$.  The choice of a particular $\rho\in\text{Hom}(\pi_{1}(N),G)$ in the conjugacy class $[\rho]$ identifies the Zariski tangent space at $\rho$ of the space $\text{Hom}(\pi_{1}(N),G)/G$  with the first cohomology group $H^{1}(N,\mathfrak{g}_{\rho})$.\\
Furthermore, the Lie algebra of the isotropy group of $\rho$
(the subgroup of $G$ fixing the representation $\rho$ under conjugation) is $H^{0}(N,\mathfrak{g}_{\rho})$.
\end{prop}
Using the definition of the R-torsion in Eq. (\ref{tor2}) above, we may define volumes on $C_{*}(N,\mathfrak{g}_{\rho})$ by using the metric on $\mathfrak{g}$.  We take $\{\sigma_{j}^{i}\otimes X_{k}\}$ to be an orthonormal basis of $C_{*}(N,\mathfrak{g}_{\rho})$, where $\sigma_{j}^{i}$ are the $i$-cells in the universal cover $\widetilde{N}$ and the $X_{k}$ are an orthonormal basis of $\mathfrak{g}$.  This volume is well defined since the adjoint representation is an orthogonal representation of $G$, and hence compatible with the equivalence relation (\ref{eqrel}).  The torsion is then an element of
\begin{equation*}
\smash[t]{\text{Tor}(C_{*}(N,\mathfrak{g}_{\rho}),\mu)\in\otimes_{2i+1}[\wedge^{\text{max}}H_{2i+1}(N,\mathfrak{g}_{\rho})]\otimes_{2i}[\wedge^{\text{max}}H_{2i}(N,\mathfrak{g}_{\rho})]^{*} }
\end{equation*}
Since $H^{i}(N,\mathfrak{g}_{\rho})\simeq H_{i}(N,\mathfrak{g}_{\rho})^{*}$, the torsion may be identified with an element
\begin{equation}\label{tors}
\smash[t]{
\text{Tor}(C_{*}(N,\mathfrak{g}_{\rho}),\mu)\in\otimes_{2i+1}[\wedge^{\text{max}}H^{2i+1}(N,\mathfrak{g}_{\rho})]^{*}\otimes_{2i}[\wedge^{\text{max}}H^{2i}(N,\mathfrak{g}_{\rho})] }
\end{equation}
The isomorphisms in Lemma \ref{conju} identify the torsion $\text{Tor}(C_{*}(N,\mathfrak{g}_{\rho}))$ with $\text{Tor}(C_{*}(N,\mathfrak{g}_{h\rho h^{-1}})$, so the torsion descends to an equivalence class $\tau(N,\rho)$ depending only on the conjugacy class $[\rho]\in\text{Hom}(\pi_{1}(N),G)/G$.\\

It is instructive to consider the case when $N$ is a genus $g$ surface $\Sigma^{g}$.
\begin{example}
From the relation (\ref{tors}) above, we see that the torsion $\tau(\Sigma^{g}; \rho)$ of a surface $\Sigma^{g}$ takes values in
\begin{equation*}
\smash[t]{
\wedge^{\rm max}H^{1}(\Sigma^{g},\mathfrak{g}_{\rho})^{*}\otimes
\wedge^{\rm max}H^{2}(\Sigma^{g},\mathfrak{g}_{\rho})\otimes\wedge^{\rm max}H^{0}(\Sigma^{g},\mathfrak{g}_{\rho}) }
\end{equation*}
By Poincare duality $H^{2}(\Sigma^{g},\mathfrak{g}_{\rho})$ is canonically dual to $H^{0}(\Sigma^{g},\mathfrak{g}_{\rho})$, so we have
\begin{equation*}
\smash[t]{\tau(\Sigma^{g}; \rho)\in\wedge^{\rm max}H^{1}(\Sigma^{g},\mathfrak{g}_{\rho})^{*} }
\end{equation*}
Observe that $H^{1}(\Sigma^{g},\mathfrak{g}_{\rho})$ is a symplectic vector space (the tangent space to the moduli space of gauge equivalence classes of
flat connections on
$\Sigma^g$) with the symplectic form
given by the cup product $H^{1}(\Sigma^{g},\mathfrak{g}_{\rho})\otimes H^{1}(\Sigma^{g},\mathfrak{g}_{\rho})\rightarrow H^{2}(\Sigma^{g},\R)\simeq \R$.  One can show that in fact the torsion may be identified with the symplectic volume on $H^{1}(\Sigma^{g},\mathfrak{g}_{\rho})$.  A rigorous proof of this is given in \cite{w}.
\end{example}
The case of interest to us here is when $N$ is a Seifert manifold and $G=U(1)$.  In this case, the torsion $\tau(N; \rho)$ takes values in
\begin{equation*}\smash[t]{
\wedge^{\text{max}}H^{1}(N,\mathfrak{g}_{\rho})^{*}\wedge^{\text{max}}H^{3}(N,\mathfrak{g}_{\rho})^{*}\otimes\wedge^{\text{max}}H^{2}(N,\mathfrak{g}_{\rho})\otimes\wedge^{\text{max}}H^{0}(N,\mathfrak{g}_{\rho})
}\end{equation*}
where by Poincar\'e duality $H^{3}(N,\mathfrak{g}_{\rho})$ is canonically dual to
\begin{equation*}\smash[t]{
H^{0}(N,\mathfrak{g}_{\rho})\simeq\R}
\end{equation*}
and $H^{1}(N,\mathfrak{g}_{\rho})$ is canonically dual to $H^{2}(N,\mathfrak{g}_{\rho})$.  Note that
\begin{equation*}
\smash[t]{H^{0}(N,\mathfrak{g}_{\rho})\simeq\R}
\end{equation*}
once we choose a basis because we are working with $U(1)$.  Thus,
\begin{equation*}
\smash[t]{
\tau(N; \rho)\in (\wedge^{\text{max}}H^{1}(N,\mathfrak{g}_{\rho})^{*})^{\otimes 2}
}\end{equation*}
or
\begin{equation*}\smash[t]{
\sqrt{\tau(N; \rho)}\in\wedge^\text{max}H^{1}(N,\mathfrak{g}_{\rho})^{*}}
\end{equation*}
 When $N$ is  two dimensional
we observed above (see \cite{w})  that the torsion can be identified with the symplectic volume on $H^{1}(\Sigma,\mathfrak{g}_{\rho})$.

We recall our previous results in equation \ref{symvol}, where we observed that for the gauge group $U(1)$ the Chern-Simons partition function  $Z_X$ was proportional to the symplectic volume. In the  case $G=U(1)$,  for a Seifert manifold $U(1)\rightarrow N\rightarrow \Sigma$, we would also like to see that $\sqrt{\tau(N)}$ is proportional  to the symplectic volume on the moduli space $U(1)^{2g}\times \Z_n$\\

Writing this more concisely, we want to see that
\begin{equation}\label{comp}\smash[t]{
\sqrt{\tau(N)}=C\cdot \omega^{g} }
\end{equation}
where $C\in\R^{*}$ is some non-zero constant.
Here $\omega\in \Omega^{2}(U(1)^{2g}\times \Z_n);\R)$ is the symplectic form on $U(1)^{2g}\times \Z_n$, i.e. the symplectic form on each of the $n$ disjoint copies of $U(1)^{2g}$.
We introduce the notation $a =  (\rho,m) \in U(1)^{2g} \times \Z_n$. Here
  $\omega_{a}(\alpha,\beta):=\int_{\Sigma}\alpha\wedge \beta$, for $\alpha,\beta\in H^{1}(N,\mathfrak{g}_{a})\simeq H^{1}(\Sigma,\mathfrak{g}_{\rho})\simeq T_{a}(U(1)^{2g}\times \Z_n).$ Here the last equality follows from Proposition \ref{jw}.
Also $\sqrt{\tau(N)}\in \Omega^{2g}(U(1)^{2g}\times \Z_n;\R)$, i.e. we define $\sqrt{\tau(N)}_{a}:=\sqrt{\tau(N; a)},\,\,\,\forall a\in U(1)^{2g}\times \Z_n$ as a section of
the top exterior power of $T^{*}(U(1)^{2g}\times \Z_n)]$.
To summarize, both the torsion and the symplectic volume are volume
elements on $U(1)^{2g} \times \Z_n$ and the Lie group structure of
$U(1)$ means the tangent bundle is trivial, and there is a natural
basis vector  given by a generator of the Lie algebra of $U(1)$. In terms of this
basis vector, we note that the definition of the torsion does not depend on
the choice of a point in $U(1)^{2g} \times \Z_n$, since the differential
of the chain complex is simply the  exterior derivative.
For nonabelian groups the differential is the twisted
differential $d_A = d + {\rm Ad} (A)$ which does depend on the choice
of a flat connection $A$. If the group is abelian $d_A$ reduces to
the exterior derivative $d$, and it does not depend on $A$.

It will be sufficient to identify $\sqrt{\tau(N)}$ and $\omega^{g}$ at a single point of the moduli space $U(1)^{2g}\times \Z_n$ since we will show that
$\sqrt{\tau(N)}$ and $\omega$ are invariant under left multiplication,  i.e. invariant under the action,
\begin{equation*} \smash[t]{
L_{a}:U(1)^{2g}\times \Z_n\rightarrow U(1)^{2g}\times \Z_n }
\end{equation*}
defined by $L_{a}:b\mapsto a\cdot b$ for $a\in U(1)^{2g}\times \Z_n$.\\

First, we show that $\sqrt{\tau(N)}$ is invariant under this action.
As discussed above, the torsion does not depend
on the choice of  point $a = [A] \in  U(1)^{2g} \times \Z_n $ corresponding
to a flat connection $A$.
 By construction,
the torsion is independent of $a$.

We identify $H^{1}(N,\R)\simeq H^{1}(N,d)$, viewed as differential $1$-forms, with $T_{a}(U(1)^{2g}\times \Z_n)$ in the following way.  The generators of $H^{1}(N,\R)$ come from generators of the cohomology for $U(1)$; one for each generating loop of the fundamental group $\pi_1(N)$.  The tangent space to the moduli space $U(1)^{2g}\times \Z_{n}$ at $a$ then has $2g$ generators $\{ \frac{\partial}{\partial \phi^{i}}_{a}\}_{i=1}^{2g}$.
Since the definition of $\tau(N)$ is independent of $a$,
where $\frac{\partial}{\partial\phi^{i}}$ is the vector field on the $i^{\text{th}}$ copy of $U(1)$ in $U(1)^{2g}\times \Z_n$, we have that
\begin{equation} \smash[t]{ 
(dL_{a})_{b} \frac{\partial}{\partial \phi^{i}}_{a}=\frac{\partial}{\partial \phi^{i}}_{a\cdot b},\,\,\forall b\in U(1)^{2g}\times \Z_{n},\,\, 1\leq i\leq 2g }
\end{equation}
We conclude that $L^{*}_a\tau(N)=\tau(N),\,\,\forall a\in U(1)^{2g}\times \Z_{n}$, i.e. $\tau(N)$ is invariant under the action of $U(1)^{2g}\times \Z_{n}$ on itself.\\

The next thing that we will show is that $\omega$, the Goldman symplectic form on the moduli space, is invariant under the action of $U(1)^{2g}\times \Z_{n}$ on itself.  This can be seen directly.  Consider
\begin{eqnarray*}
(L_{a}^{*}\omega)_{a}(\frac{\partial}{\partial \phi^{i}},\frac{\partial}
{\partial \theta^{j}})&=&\omega_{a\cdot a}((dL_a)_{a}(\frac{\partial}{\partial \phi^{i}}),(dL_a)_{a}(\frac{\partial}{\partial \phi^j}))\\
                                                         &=&\omega_{a\cdot a}(\frac{\partial}{\partial \phi^{i}},\frac{\partial}{\partial \phi^j})\\
                                                         &=&\omega_{a}(\frac{\partial}{\partial \phi^{i}},\frac{\partial}{\partial \phi^{j}}.)
\end{eqnarray*}

Thus,
\begin{equation*}\smash[t]{
L_{a}^{*}\omega=\omega,\,\,\forall a\in U(1)^{2g}\times \Z_{n} }
\end{equation*}
i.e. $\omega$ is invariant under the action of $U(1)^{2g}\times \Z_{n}$ on itself.

Now we can prove our original claim.  Let $e$ denote the identity element of $U(1)^{2g}\times \Z_{n}$.  Then at the point $e$, $\sqrt{\tau(N)}$ and $\omega^{g}$ must agree up to a non-zero multiplicative constant:
\begin{equation*}\smash[t]{
\sqrt{\tau(N)}|_{e}=C\cdot \omega^{g}|_{e}}
\end{equation*}
for some $C\in\R^{*}$.  By left invariance, we therefore have:
\begin{equation*}\smash[t]{
\sqrt{\tau(N)}|_{a}=C\cdot \omega^{g}|_{a},\,\,\,\forall a\in U(1)^{2g}\times \Z_n }
\end{equation*}
Thus,
\begin{equation*} 
\sqrt{\tau(N)}=C\cdot \omega^{g}
\end{equation*}



\addcontentsline{toc}{chapter}{Bibliography}


\end{document}